\documentclass[11pt,english]{smfart}
\usepackage{amssymb}
\usepackage{latexsym}

\def\di{\displaystyle}

\font\bbb=msbm10 at 12pt
\def\rR{\hbox{\bbb R}}

\def\cC{\hbox{\bbb C}}

\newtheorem{thm}{Theorem}[section]
\newtheorem{defi}{Definition}[section]
\newtheorem{lem}{Lemma}[section]

\newtheorem{rema}{Remark}[section]

\title{Non differentiable variational principles}
\alttitle{Principes variationnels non diff\'erentiables}
\author{Jacky CRESSON}
\address{Universit\'e de Franche-Comt\'e, Equipe de Math\'ematiques de Besan\c{c}on,
CNRS-UMR 6623, Th\'eorie des nombres et alg\`ebre, 16 route de
gray, 25030 Besan\c{c}on cedex, France.}
\email{cresson@math.univ-fcomte.fr}

\begin{document}

\maketitle
\baselineskip 6mm

\begin{abstract}
We developp a calculus of variations for functionals which are defined on a set of non differentiable
curves. We first extend the classical differential calculus in a quantum calculus, which allows us to
define a complex operator, called the scale derivative, which is the non differentiable analogue of
the classical derivative. We then define the notion of extremals for our functionals and obtain a
characterization in term of a generalized Euler-Lagrange
equation. We finally prove that solutions of the Schr\"odinger equation can be
obtained as extremals of a non differentiable variational principle, leading to an extended
Hamilton's principle of least action for quantum mechanics. We compare this approach with the
scale relativity theory of Nottale, which assumes a fractal structure of space-time.

\end{abstract}

\begin{altabstract}
Nous d\'eveloppons un calcul des variations pour des fonctionnelles d\'efinies sur un ensemble de
courbes non diff\'erentiables. Pour cela, nous \'etendons le calcul diff\'erentiel classique, en
calcul appel\'e {\it calcul quantique}, qui nous permet de d\'efinir un op\'erateur \`a valeur complexes,
appel\'e {\it d\'eriv\'ee d'\'echelle}, qui est l'analogue non diff\'erentiable de la d\'eriv\'ee
usuelle. On d\'efinit alors la notion d'extremale pour ces fonctionnelles pour lesquelles nous obtenons
une caract\'erisation via une \'equation d'Euler-Lagrange g\'en\'eralis\'ee. On prouve enfin que
les solutions de l'\'equation de Schr\"odinger peuvent s'obtenir comme solution d'un probl\`eme
variationnel non diff\'erentiable, \'etendant ainsi le principe de moindre action de Hamilton au
cadre de la m\'ecanique quantique. On discute enfin la connexion entre ce travail et la
th\'eorie de la relativit\'e d'\'echelle d\'evelopp\'ee par Nottale, et qui suppose une structure
fractale de l'espace-temps.
\end{altabstract}

\tableofcontents

\section{Introduction}

Lagrangian mechanics describes motion of mechanical systems using differentiable manifolds. Motions
of Lagrangian systems are extremals of a variational principle called ``Hamilton's principle
of least action" (see \cite{ar},p.55).\\

However, some important physical systems can't be put in such a framework. For example, generic
trajectories of quantum mechanics are not differentiable curves \cite{fh}, such that a classical
Lagrangian formalism is not possible (see however \cite{fey}).\\

In this article we extend the calculus of variations in order to cover sets of non differentiable
curves. We first define a quantum calculus allowing us to analyze non differentiable functions
by means of a complex operator, which generalizes the classical derivative. We then introduce
functionals on H\"olderian curves and study the analogue of extremals for these objects. We
prove that extremals curves of our functionals are solutions of a generalized Euler-Lagrange
equation, which looks like the one obtain by Nottale \cite{no} in the context of the scale relativity
theory. We then prove that the Schr\"odinger equation can be obtain as extremals of a non differentiable
variational problem.\\

The non differentiable calculus of variations gives a rigorous basis to the scale relativity principle
developped by Nottale \cite{no} in order to recover quantum mechanics by keeping out the
differentiability assumption of the space-time.

\section{Quantum calculus}

In this section we define the quantum calculus, which extends the classical differential calculus
to non differentiable functions. We refer to \cite{bc2} and \cite{cr2} for analogous ideas and the
underlying physical framework leading to this extension.

\subsection{Basic definitions}

We denote by $C^0$ the set of continuous real valued functions defined on $\rR$.

\begin{defi}
Let $f\in C^0$. For all $\epsilon >0$, we call $\epsilon$ left and right quantum derivatives the quantities
\begin{equation}
\Delta^{\sigma}_{\epsilon} f(t)=\sigma \di {f(t+\sigma \epsilon )-f(t) \over \epsilon} ,\ \sigma =\pm .
\end{equation}
\end{defi}

The $\epsilon$ left and right quantum derivatives of a continuous function correspond to the classical derivatives
of the left and right $\epsilon$-mean function defined by
\begin{equation}
f^{\sigma}_{\epsilon } (t)=\di {\sigma \over \epsilon}\di\int_t^{t+\sigma \epsilon}
f(s) ds,\ \sigma=\pm .
\end{equation}

Using $\epsilon$ left and right derivatives, we can define an operator which generalize
the classical derivative.

\begin{defi}
Let $f\in C^0$. For all $\epsilon >0$, the $\epsilon$ scale derivative of $f$ at point $t$ is the quantity denoted by
$\Box_{\epsilon} f /\Box t$, and defined by
\begin{equation}
\di {\Box_{\epsilon} f\over \Box t} (t)=\di (\Delta^+_{\epsilon} f(t)+\Delta^-_{\epsilon} f(t))-i(
\Delta^+_{\epsilon} f(t)-\Delta^-_{\epsilon} f(t)) .
\end{equation}
\end{defi}

If $f$ is differentiable, we can take the limit of the scale derivative when $\epsilon$ goes to
zero. We then obtain the classical derivative of $f$, $f'$.\\

In the following, we will frequently denote $\Box_{\epsilon} x$ for $\Box_{\epsilon} x /\Box t$.\\

We also need to extend the scale derivative to complex valued functions.

\begin{defi}
Let $f$ be a continuous complex valued function. For all $\epsilon >0$, the $\epsilon$ scale
derivative of $f$, denoted by $\Box_{\epsilon} f /\Box t$ is defined by
\begin{equation}
\di {\Box_{\epsilon} f \over \Box t} (t) =\di {\Box_{\epsilon} \mbox{\rm Re}(f) \over \Box t} +
i {\Box_{\epsilon} \mbox{\rm Im}(f) \over \Box t} ,
\end{equation}
where $\mbox{\rm Re}(f)$ and $\mbox{\rm Im}(f)$ denote the real and imaginary part of $f$.
\end{defi}

This extension of the scale derivative in order to cover complex valued functions is far from
being trivial. Indeed, it mixes complex terms in a complex operator.

\subsection{Basic formulas}

For all $\epsilon >0$ the scale derivative is not a derivation\footnote{We recall that a derivation on
an abtract algebra $A$ is a linear application $D:A\rightarrow A$ such that $D(xy)=D(x).y+x.D(y)$ for all
$x,y\in A$.} on the set of continuous functions\footnote{A classical result says that there
exists no derivations on the set of continuous functions except the trivial one, define by
$D(f)=0$ for all $f\in C^0$.} Indeed, we have  :

\begin{thm}
Let $f$ and $g$ be two functions of $C^0$. For all $\epsilon >0$ we have
\begin{equation}
\label{genleib}
\Box_{\epsilon} (fg)=\Box_{\epsilon} f .g +f.\Box_{\epsilon} g +\epsilon i\left [ \Box_{\epsilon} f
\boxminus_{\epsilon} g -\boxminus_{\epsilon} f \Box_{\epsilon} g -
\Box_{\epsilon} f \Box_{\epsilon} g -\boxminus_{\epsilon} f \boxminus_{\epsilon} g  \right ],
\end{equation}
where $\boxminus f$ is the complex conjugate of $\Box f$.
\end{thm}

Of course, when we restrict our attention to differntiable functions, taking the limit of (\ref{genleib})
when $\epsilon$ goes to zero, we obtain the classical Leibniz rule $(fg)'=f' .g+f.g'$.\\

\begin{proof}
Formula (\ref{genleib}) follows from easy calculations. In particular, we use the fact that
\begin{equation}
\Delta_{\sigma}^{\epsilon} (fg)=\Delta_{\sigma}^{\epsilon} f .g +f.\Delta_{\sigma}^{\epsilon} g +\sigma
\epsilon \Delta_{\sigma}^{\epsilon} f .\Delta_{\sigma}^{\epsilon} g  ,\ \sigma =\pm ,
\end{equation}
which is a standard result of the calculus of finite differences (see \cite{mt}).\\

As a consequence, we have
\begin{equation}
\Box_{\epsilon} (fg)=\Box_{\epsilon} f .g +f.\Box_{\epsilon} g +\epsilon \left [ (\Delta_{\epsilon}^+ f
\Delta_{\epsilon}^- g - \Delta_{\epsilon}^- f \Delta_{\epsilon}^- g) -i
(\Delta_{\epsilon}^+ f \Delta_{\epsilon}^- g + \Delta_{\epsilon}^- f \Delta_{\epsilon}^- g)
\right ]
.
\end{equation}
Moreover, we have the following formula :
\begin{equation}
\Box_{\epsilon } f \Box_{\epsilon} g =\di {1\over 2} \left [ (\Delta_{\epsilon}^+ f \Delta_{\epsilon}^- g +\Delta_{\epsilon}^- f
\Delta_{\epsilon}^+ g ) -i ( \Delta_{\epsilon}^+ f \Delta_{\epsilon}^+ g -\Delta_{\epsilon}^- f
\Delta_{\epsilon}^- g )\right ] ,
\end{equation}
and
\begin{equation}
\Box_{\epsilon } f \boxminus_{\epsilon} g =\di {1\over 2} \left [ ( \Delta_{\epsilon}^+ f \Delta_{\epsilon}^+ g
+\Delta_{\epsilon}^- f \Delta_{\epsilon}^- g )
-i(\Delta_{\epsilon}^+ f \Delta_{\epsilon}^- g  -\Delta_{\epsilon}^- f
\Delta_{\epsilon}^+ g )
\right ] .
\end{equation}
We then obtain
\begin{equation}
\left .
\begin{array}{lll}
\Delta_{\epsilon}^+ f \Delta_{\epsilon}^+ g
+\Delta_{\epsilon}^- f \Delta_{\epsilon}^- g  & = & \Box_{\epsilon } f \boxminus_{\epsilon} g +
\boxminus_{\epsilon} f \Box_{\epsilon} g ,\\
-i(\Delta_{\epsilon}^+ f \Delta_{\epsilon}^+ g
-\Delta_{\epsilon}^- f \Delta_{\epsilon}^- g)  & = & \Box_{\epsilon } f \Box_{\epsilon} g -
\boxminus_{\epsilon} f \boxminus_{\epsilon} g .
\end{array}
\right .
\end{equation}
We deduce then the following equality
\begin{equation}
\left .
\begin{array}{l}
(\Delta_{\epsilon}^+ f
\Delta_{\epsilon}^- g - \Delta_{\epsilon}^- f \Delta_{\epsilon}^- g) -i
(\Delta_{\epsilon}^+ f \Delta_{\epsilon}^- g + \Delta_{\epsilon}^- f \Delta_{\epsilon}^- g)
\\
=i (\Box_{\epsilon } f \Box_{\epsilon} g -
\boxminus_{\epsilon} f \boxminus_{\epsilon} g)
-i (\Box_{\epsilon } f \boxminus_{\epsilon} g +
\boxminus_{\epsilon} f \Box_{\epsilon} g) .
\end{array}
\right .
\end{equation}
This concludes the proof.
\end{proof}

We have the following integral formula :
\begin{equation}
\label{formint}
\di\int_a^b \Box_{\epsilon} f (t) dt =\left . \di {1\over 2} \left [ \left ( f^+_{\epsilon} (t)+f^-_{\epsilon} (t) \right )
-i \left ( f^+_{\epsilon} (t)-f^-_{\epsilon} (t) \right  ) \right ] \right |_a^b  .
\end{equation}
When $\epsilon$ goes to zero, we deduce
\begin{equation}
\lim_{\epsilon \rightarrow 0} \di\int_a^b \Box_{\epsilon} f(t) dt = \left . f(t)\right |_a^b .
\end{equation}

\subsection{H\"olderian functions}

In the following, we consider a particular class of non differentiable functions called
{\it H\"olderian functions} \cite{tri}.

\begin{defi}
A continuous real valued function $f$ is H\"olderian of H\"older exponent $\alpha$, $0<alpha <1$, if for
all $\epsilon >0$, and all $t,t'\in \rR$ such that $\mid t-t'\mid \leq \epsilon$, there exists a
constant $c$ such that
\begin{equation}
\label{hold}
\mid f(t) -f(t')\mid \leq c \epsilon^{\alpha} .
\end{equation}
\end{defi}

In the following, we denote by $H^{\alpha}$ the set of continuous functions which are H\"olderian
of H\"older exponent $\alpha$. Moreover, we say that a complex valued function $y(t)$ belongs
to $H^{\alpha}$ if its real and imaginary part belong to $H^{\alpha}$.\\

We then have the following lemma :

\begin{lem}
\label{regsd}
If $x\in H^{\alpha}$ then $\Box_{\epsilon} x \in H^{\alpha}$ for all $\epsilon >0$.
\end{lem}

This follows from the definition of $\Box_{\epsilon} x(t)$ and simple calculations.

\subsection{A technical result}

We derive a technical result about the scale derivative, which will be used in the last section.

\begin{thm}
\label{tech}
Let $f(x,t)$ be a $C^{n+1}$ function and $x(t)\in H^{1/n}$, $n\geq 1$. For all $\epsilon >0$ sufficiently small,
we have
\begin{equation}
\di {\Box_{\epsilon} f (x(t),t) \over \Box t} =
\di {\partial f\over \partial t} +
\di \sum_{j=1}^n \di {1\over j!} \di {\partial^j f \over \partial x^j } (x(t),t)
\epsilon^{j-1} a_{\epsilon ,j} (t) +o (\epsilon^{1/n} ),
\end{equation}
where
\begin{equation}
a_{\epsilon ,j} (t) =
\di {1\over 2}
\left [
\left (
(\Delta_+^{\epsilon} x )^j -(-1)^j (\Delta_-^{\epsilon} x )^j
\right )
-i
\left (
(\Delta_+^{\epsilon} x )^j +(-1)^j (\Delta_-^{\epsilon} x )^j
\right )
\right ]
.
\end{equation}
\end{thm}

The proof follows easily from the following lemma :

\begin{lem}
Let $f(x,t)$ be a real valued function of class $C^{n+1}$, $n\geq 1$, and $x(t)\in H^{1/n}$.
For all $\epsilon >0$ sufficiently small, the right and left quantum derivatives of $f(x(t),t)$ are
given by
\begin{equation}
\Delta_{\sigma}^{\epsilon} f(x(t),t) =\di {\partial f \over \partial t} (x(t),t) +\sigma
\di\sum_{i=1}^n
\di {1\over i!} \di {\partial^i f \over \partial x^i} (x(t),t) \epsilon^{-1} (\sigma
\epsilon \Delta_{\sigma}^{\epsilon} x(t))^i
+o(\epsilon^{1/n} ) ,
\end{equation}
for $\sigma =\pm$.
\end{lem}

\begin{proof}
This follows from easy computations. First, we remark that, as $x(t)\in H^{1/n}$,
we have $\mid \epsilon \Delta_{\sigma}^{\epsilon} X(t) \mid =o(\epsilon^{1/n} )$. Moreover,
$$
f(x(t+\epsilon ) ,t+\epsilon ) =f(x(t)+\epsilon \Delta_+^{\epsilon} x(t), t+\epsilon ) .
$$
By the previous remark, and the fact that $f$ is of order $C^{n+1}$, we can make a Taylor
expansion up to order $n$ with a controled remainder.
$$
\left .
\begin{array}{lll}
f(x(t+\epsilon) ,t+\epsilon )& = & f(x(t),t)+\di\sum_{k=1}^{n} \di {1\over k!} \di \sum_{i+j=k}
(\epsilon \Delta_+^{\epsilon} x(t))^i \epsilon^j \di {\partial^k f \over \partial^i x \partial^j t}
(x(t),t) \\
 & & +o((\epsilon \Delta_+^{\epsilon} x(t))^{n+1} ) .
\end{array}
\right .
$$
As a consequence, we have
$$
\epsilon \Delta_+^{\epsilon} f(x(t),t) =\di\sum_{k=1}^{n} \di {1\over k!} \di \sum_{i+j=k}
(\epsilon \Delta_+^{\epsilon} x(t))^i \epsilon^j \di {\partial^k f \over \partial^i x \partial^j t}
(x(t),t)
+o((\epsilon \Delta_+^{\epsilon} x(t))^{n+1} ) .
$$
By selecting terms of order less or equal to one in $\epsilon$ in the right of this equation, we
obtain
$$
\left .
\begin{array}{lll}
\epsilon \Delta_+^{\epsilon} f(x(t),t) & = & \epsilon \left [
\di {\partial f \over \partial t} (x(t),t) +\di\sum_{i=1}^n
\di {1\over i!} \di {\partial^i f \over \partial x^i} (x(t),t) \epsilon^{-1} (\epsilon
\Delta_+^{\epsilon} x(t))^i
\right ] \\
 & & +o(\epsilon^2 \Delta_+^{\epsilon} x(t) ) .
\end{array}
\right .
$$
Dividing by $\epsilon$, we obtain the lemma.
\end{proof}

\section{Non differentiable calculus of variations}

\subsection{Functionals}

The classical calculus of variations is concerned with the extremals of functions whose domain is an
infinite-dimensional space : the space of curves, which is usually the set of {\it differentiable}
curves. We look for an analogous theory on the set of non differentiable curves.\\

In all the text, $\alpha$ is a real number satisfying
$$0<\alpha <1 ,$$
and $\epsilon$ is a parameter, which is assumed to be sufficiently small, i.e.
$$0<\epsilon <<1,$$
without precising its exact smallness.

We denote by $C_{\epsilon}^{\alpha} (a,b)$ the set of curves in the plane of the form
\begin{equation}
\gamma =\{ (t,x(t)), x\in H^{\alpha} , a-\epsilon \leq t\leq b+\epsilon \}   .
\end{equation}

\begin{rema}
i. In the following, we will simply write $C^{\alpha} (a,b)$ for $C_{\epsilon}^{\alpha} (a,b)$.\\

ii. We must take $a-\epsilon \leq t\leq b+\epsilon$ in order to avoid problems with the definition
of the scale derivative on the extremal points of the interval $[a,b]$.
\end{rema}

A {\it functional} $\Phi$ is a map $\Phi :C^{\alpha} (a,b) \rightarrow \cC$.

\begin{rema}
In classical mechanics, one usually consider real valued functionals instead of complex one.
\end{rema}

We will restrict our attention to the following class of functionals :

\begin{defi}
Let $L :\rR \times \cC \times \rR \rightarrow \cC$ be a differentiable function of three
variables $(x,v,t)$. For all $\epsilon >0$, a functional $\Phi_{\epsilon} :C^{\alpha} (a,b)\rightarrow \cC$ is defined by
\begin{equation}
\label{funct}
\Phi_{\epsilon} (\gamma )=\int_a^b L(x(t),\Box_{\epsilon} x(t) ,t) dt ,
\end{equation}
for all $\gamma \in C^{\alpha} (a,b)$.
\end{defi}

Of course, when we consider differentiable curves, we can take the limit of (\ref{funct}) when
$\epsilon$ goes to zero, and we obtain the classical functional (see \cite{ar},p.56) :
\begin{equation}
\Phi (\gamma )=\int_a^b L(x(t),\dot{x}(t),t) dt ,
\end{equation}
where $\dot{x}=dx/dt$.

\subsection{Variations}

We first define {\it variations} of curves.

\begin{defi}
Let $\gamma \in C^{\alpha} (a,b)$. A variation $\gamma'$ of $\gamma$ is a curve
\begin{equation}
\gamma' =\{ (t,x(t)+h(t)), x\in H^{\alpha} , h\in H^{\beta} ,
\beta \geq \alpha 1_{[1/2 ,1]} +(1-\alpha )1_{]0,1/2[} , h(a)=h(b)=0 \} .
\end{equation}
We denote this curve by $\gamma' =\gamma +h$.
\end{defi}

As in the usual case, we look for paths of a given regularity class with prescribed end points. The
condition $\beta \geq \alpha 1_{[1/2 ,1]} +(1-\alpha )1_{]0,1/2[}$ for the variation is
a technical assumption, which will be used in the derivation of the  non differentiable analogue of
the Euler-Lagrange equation (see $\S$.\ref{eulerlag}). The minimal condition on $\beta$ for which the
problem of variations makes sense is $\beta \geq \alpha$, in order to ensure that $\gamma +h$ is again
in $C^{\alpha} (a,b)$.\\

In the following, we always consider variations of a given curve $\gamma$ of the form
$\gamma_{\mu} =\gamma +\mu h$, where $\mu$ is a real parameter.

\begin{defi}
A functional $\Phi$ is called differentiable on $C^{\alpha} (a,b)$ if for all variations
$h\in C^{\beta} (a,b)$, we have
\begin{equation}
\Phi (\gamma +h )-\Phi (\gamma)=F (\gamma ,h)+R(\gamma ,h) ,
\end{equation}
where $F$ depends linearly on $h$, i.e. $F(h_1 +h_2 )=F(h_1 )+F(h_2)$ and $F(ch)=cF(h)$, and
$R(\gamma ,h)=O(h^2 )$, i.e. for $\mid h\mid <\mu$ and $\mid \Box_{\epsilon} h\mid <\mu$, we
have $\mid R\mid <C\mu^2$.

The functional $F$ is called the differential of $\Phi$.
\end{defi}

In the case of functionals of the form (\ref{funct}), we have :

\begin{thm}
\label{funcder}
For all $\epsilon >0$, the functional $\Phi_{\epsilon} (\gamma )$ defined by (\ref{funct}) is differentiable,
and its derivative is given by the formula
\begin{equation}
\label{eq1}
\left .
\begin{array}{lll}
F_{\epsilon}^{\gamma} (h) & = & \di\int_a^b
\left [ \di {\partial L \over \partial x} (x(t), {\Box_{\epsilon} x \over \Box t} ,t)
-\di {\Box_{\epsilon} \over \Box t} \left (
\di {\partial L \over \partial \Box_{\epsilon} x} (x(t),{\Box_{\epsilon} x \over \Box t} ,t)
\right ) \right ] h (t) dt \\
 & &  + \di\int_a^b \di {\Box_{\epsilon} \over \Box t}  \left (
 \di {\partial L\over \partial \Box_{\epsilon} x} h (t) \right ) dt
+iR_{\epsilon}^{\gamma} (h) ,
\end{array}
\right .
\end{equation}
with
\begin{equation}
R_{\epsilon}^{\gamma} (h) =
\epsilon \di\int_a^b \di \left [
\Box_{\epsilon} f_{\epsilon} (t) \Box_{\epsilon} h(t)
-
\boxminus_{\epsilon} f_{\epsilon} (t) \Box_{\epsilon} h(t)
-
\Box_{\epsilon} f_{\epsilon} (t) \boxminus_{\epsilon} h(t)
-
\boxminus_{\epsilon} f_{\epsilon} (t) \boxminus_{\epsilon} h(t)
\right ] dt ,
\end{equation}
where
\begin{equation}
f_{\epsilon} (t) =\di {\partial L\over \partial \Box_{\epsilon} x} (x(t),\Box_{\epsilon} x(t),t).
\end{equation}
\end{thm}

\begin{proof}
We have
$$
\left .
\begin{array}{l}
\Phi_{\epsilon} (\gamma +h) -\Phi (\gamma ) = \\
\di\int_a^b \left [ L ( x(t)+ h (t) ,\Box_{\epsilon} x (t) + \Box_{\epsilon} h (t)
         ,t ) -L(x(t),\Box_{\epsilon} x(t),t) \right ] dt  ,\\
 =  \di\int_a^b \left [ \di {\partial L\over \partial x } (x(t),\Box_{\epsilon} x(t),t) h (t)
 + \di {\partial L\over \partial \Box_{\epsilon} x } (x(t),\Box_{\epsilon} x(t),t) \Box_{\epsilon}
 h (t ) \right ] dt + O(h^2 ) ,\\
   =  F_{\epsilon}^{\gamma} (h)+R(h),
\end{array}
\right .
$$
where
$$
F_{\epsilon}^{\gamma} (h) =\di\int_a^b \left [ \di {\partial L\over \partial x } (x(t),\Box_{\epsilon} x(t),t) h (t)
 + \di {\partial L\over \partial \Box_{\epsilon} x } (x(t),\Box_{\epsilon} x(t),t) \Box_{\epsilon}
 h (t ) \right ] dt,
$$
and $R(h)=O(h^2 )$.\\

Using (\ref{genleib}), we deduce :
\begin{equation}
\label{eq1}
\left .
\begin{array}{lll}
F_{\epsilon}^{\gamma} (h) & = & \di\int_a^b
\left [ \di {\partial L \over \partial x} (x(t), {\Box_{\epsilon} x \over \Box t} ,t)
-\di {\Box_{\epsilon} \over \Box t} \left (
\di {\partial L \over \partial \Box_{\epsilon} x} (x(t),{\Box_{\epsilon} x \over \Box t} ,t)
\right ) \right ] h (t) dt \\
 & &  + \di\int_a^b \di {\Box_{\epsilon} \over \Box t}  \left (
 \di {\partial L\over \partial \Box_{\epsilon} x} h (t) \right ) dt\\
 & & +i\epsilon \di\int_a^b
\left [
\Box_{\epsilon} f_{\epsilon} (t) \Box_{\epsilon} h(t)
-
\boxminus_{\epsilon} f_{\epsilon} (t) \Box_{\epsilon} h(t)
-
\Box_{\epsilon} f_{\epsilon} (t) \boxminus_{\epsilon} h(t)
-
\boxminus_{\epsilon} f_{\epsilon} (t) \boxminus_{\epsilon} h(t)
\right ] dt ,
\end{array}
\right .
\end{equation}
with $f_{\epsilon} (t) =\di {\partial L\over \partial \Box_{\epsilon} x} (x(t),\Box_{\epsilon} x(t),t)$.
This concludes the proof.
\end{proof}

\subsection{Extremal curves and Euler-Lagrange equation}
\label{eulerlag}

The functional derivative of $\Phi_{\epsilon}$ mix terms which are either divergent when $\epsilon$
goes to zero, or tending toward $0$ with $\epsilon$. In order to simplify our problem and
to take into account only dominant terms in $\epsilon$, we introduce the following operator :

\begin{defi}
Let $a_p (\epsilon )$ be a real or complex valued function, with parameters $p$. We denote  by
$[.]_{\epsilon}$ the linear operator defined by :\\

i. $a_p (\epsilon )-[a_p (\epsilon )]_{\epsilon} \longrightarrow_{\epsilon \rightarrow 0} 0$,

ii. $[a_p (\epsilon )]_{\epsilon} =0$ if $\lim_{\epsilon \rightarrow 0} a_p (\epsilon )=0$.\\

\noindent The quantity $[a_p (\epsilon )]_{\epsilon}$ is called the $\epsilon$-dominant part of $a_p
(\epsilon )$.
\end{defi}

For example, if $a(\epsilon )=\epsilon^{-1/2} +2\epsilon +2$, then $[a(\epsilon )]_{\epsilon } =
\epsilon^{-1/2} +2$.\\

We deduce the following properties :

\begin{lem}
The $\epsilon$-dominant part is unique.
\end{lem}

\begin{proof}
This comes from the relation $[[.]_{\epsilon} ]_{\epsilon} =[.]_{\epsilon}$. Indeed, by definition we
have $a_p (\epsilon )=[a_p (\epsilon )]_{\epsilon } +r(\epsilon )$ with $\lim_{\epsilon
\rightarrow 0} r(\epsilon )=0$. Applying $[.]_{\epsilon}$ directly on this expression, we obtain
$[a_p (\epsilon )]_{\epsilon} =[[a_p (\epsilon )]_{\epsilon} ]_{\epsilon}$ using ii.
\end{proof}

\begin{rema}
Unicity comes from condition ii. Indeed, if we cancel this condition, we can obtain many different
quantities satisfying i. For example, if $a (\epsilon )=\alpha  \epsilon^{-1/2} +\epsilon +2$, then
without ii), we have the choice between $[a(\epsilon )]_{\epsilon} =\epsilon^{-1/2} +2+\epsilon$ and
$[a(\epsilon )]_{\epsilon} =\epsilon^{-1/2} +2$.

ii. This operator can be used in the definition of left and right quantum operators by considering
$\delta_{\epsilon}^{\sigma} x(t)=[\Delta_{\epsilon}^{\sigma} x(t)]_{\epsilon}$, $\sigma =\pm$. However,
using such kind of operators lead to many difficulties from the algebraic
point of view, in particular with the derivation of the analogue of the Leibniz rule.
\end{rema}

We now introduce the non differentiable analogue of the notion of {\it extremals} curves in the
classical case (see \cite{ar},p.57 ).

\begin{defi}
Let $0<\alpha \leq 1$. An extremal curve of the functional (\ref{funct}) on the space of
curves of class $C^{\beta} (a,b)$, $\beta \geq \alpha 1_{[1/2 ,1]} +(1-\alpha )1_{]0,1/2[}$,
is a curve $\gamma \in C^{\alpha} (a,b)$ satisfying
\begin{equation}
\label{prinmac}
[F_{\epsilon}^{\gamma} (h) ]_{\epsilon} =0 ,
\end{equation}
for all $\epsilon >0$ and all $h \in C^{\beta} (a,b)$.
\end{defi}

The following theorem gives the analogue of the Euler-Lagrange equations for extremals of our
functionals.

\begin{thm}
\label{eulerlagequa}
We assume that the function $L$ defining the functional (\ref{funct}) satisfies
\begin{equation}
\label{tech1}
\parallel D (\partial L /\partial v )\parallel \leq C ,
\end{equation}
where $C$ is a constant, $D$ denotes the differential, and $\parallel .\parallel$ is the classical
norm on matrices.

The curve $\gamma :\ x=x(t)$ is an extremal curve of the functional (\ref{funct}) on the space
of curves of class $C^{\beta} (a,d)$,
\begin{equation}
\label{tech3}
\beta \geq \alpha 1_{[1/2 ,1]} +(1-\alpha )1_{]0,1/2[} ,
\end{equation}
if and only if it satisfies the following generalized Euler-Lagrange equation
\begin{equation}
\label{stat}
\left [ \di {\partial L \over \partial x} (x(t), {\Box_{\epsilon} x \over \Box t} ,t)
-\di {\Box_{\epsilon} \over \Box t} \left (
\di {\partial L \over \partial \Box_{\epsilon} x} (x(t),{\Box_{\epsilon} x \over \Box t} ,t)
\right ) \right ]_{\epsilon} =0 ,
\end{equation}
for $\epsilon >0$.
\end{thm}

\begin{rema}
Our Euler-Lagrange equation (\ref{stat}) looks like the one obtain by Nottale \cite{no} in
the context of the {\it scale relativity theory} (see $\S$.\ref{srt}).
\end{rema}

\begin{proof}
The proof follow the classical derivation of Euler-Lagrange equation (see for example
\cite{spi},p.432-434). By theorem \ref{funcder}, we have
$$
\left .
\begin{array}{lll}
F_{\epsilon}^{\gamma} (h) & = & \di\int_a^b
\left [ \di {\partial L \over \partial x} (x(t), {\Box_{\epsilon} x \over \Box t} ,t)
-\di {\Box_{\epsilon} \over \Box t} \left (
\di {\partial L \over \partial \Box_{\epsilon} x} (x(t),{\Box_{\epsilon} x \over \Box t} ,t)
\right ) \right ] h (t) dt \\
 & &  + \di\int_a^b \di {\Box_{\epsilon} \over \Box t}  \left (
 \di {\partial L\over \partial \Box_{\epsilon} x} h (t) \right ) dt\\
 & & +i\epsilon \di\int_a^b
\left [
\Box_{\epsilon} f_{\epsilon} (t) \Box_{\epsilon} h(t)
-
\boxminus_{\epsilon} f_{\epsilon} (t) \Box_{\epsilon} h(t)
-
\Box_{\epsilon} f_{\epsilon} (t) \boxminus_{\epsilon} h(t)
-
\boxminus_{\epsilon} f_{\epsilon} (t) \boxminus_{\epsilon} h(t)
\right ] dt ,
\end{array}
\right .
$$
with $f_{\epsilon} (t) =\di {\partial L\over \partial \Box_{\epsilon} x} (x(t),\Box_{\epsilon} x(t),t)$.\\

In order to conclude, we need the following lemma :

\begin{lem}
\label{simpli}
Let $0<\epsilon$, $a,b\in \rR$, $h \in H^{\beta}$, $\beta \geq
\alpha 1_{[1/2 ,1]} +(1-\alpha )1_{]0,1/2[}$, such that $h (a)=h  (b)=0$,
and $f_{\epsilon} :\rR \rightarrow \cC$ such that
\begin{equation}
\label{tech2}
\sup_{s\in \{ t,t+\sigma \epsilon \}} \mid f_{\epsilon} (s)\mid \leq C \epsilon^{\alpha -1} ,
\end{equation}
for all $t\in [a,b]$.

Then, we have
\begin{equation}
\label{tech5}
\int_a^b \di {\Box_{\epsilon} \over \Box t} (f_{\epsilon} (t) h (t) ) dt =O(\epsilon^{\alpha +\beta -1  }),\
\mbox{\rm and}\
\epsilon \int_a^b \di Op_{\epsilon} (f_{\epsilon} )\di Op'_{\epsilon} (h)
dt =O(\epsilon^{\alpha +\beta } ) .
\end{equation}
where $Op_{\epsilon} $ and $Op'_{\epsilon}$ are either $\Box_{\epsilon}$ or $\boxminus_{\epsilon}$.
\end{lem}

The proof is given in the next section.\\

Using condition (\ref{tech1}), we obtain
$$
\sup_{s\in \{t,t+\sigma \epsilon\}} \mid \partial L/\partial \Box_{\epsilon} x \mid
\leq C' \epsilon^{\alpha -1} ,
$$
as $\sup_{s\in \{t,t+\sigma \epsilon\}} \left [ \max (\mid x(s)\mid ,\mid \Box_{\epsilon} x(s)\mid ,
\mid s\mid ) \right ] \leq
C" \epsilon^{\alpha -1}$.\\

Using lemma \ref{simpli} with $f_{\epsilon} (s)=(\partial L /\partial \Box_{\epsilon} t )
(x(s),\Box_{\epsilon} (s),s)$, and condition (\ref{tech3}), we deduce that
$$
\lim_{\epsilon \rightarrow 0} \di\int_a^b \di {\Box_{\epsilon} \over \Box t}  \left (
 \di {\partial L\over \partial \Box_{\epsilon} x} h (t) \right ) dt =0\
 \mbox{\rm and}\
\lim_{\epsilon \rightarrow 0} \epsilon \di\int_a^b \di Op_{\epsilon} \left (
\di {\partial L\over \partial \Box_{\epsilon} x} \right ) Op'_{\epsilon} (h(t)) dt =0 ,
$$
for $Op_{\epsilon}$ and $Op'_{\epsilon}$ which are either $\Box_{\epsilon}$ or $\boxminus_{\epsilon}$.\\

Hence, applying the operator $[.]_{\epsilon}$, we obtain
$$
\left .
\begin{array}{lll}
[F_{\epsilon}^{\gamma}(h)]_{\epsilon} & = & \left [ \di\int_a^b
\left [ \di {\partial L \over \partial  x} (x(t), {\Box_{\epsilon} x \over \Box t} ,t)
-\di {\Box_{\epsilon} \over \Box t} \left (
\di {\partial L \over \partial \Box_{\epsilon} x} (x(t),{\Box_{\epsilon} x \over \Box t} ,t)
\right ) \right ] h (t) dt \right ]_{\epsilon} ,\\
  & = & \di\int_a^b
\left [ \di {\partial L \over \partial x} (x(t), {\Box_{\epsilon} x \over \Box t} ,t)
-\di {\Box_{\epsilon} \over \Box t} \left (
\di {\partial L \over \partial \Box_{\epsilon} x} (x(t),{\Box_{\epsilon} x \over \Box t} ,t)
\right ) \right ]_{\epsilon} h (t) dt .
\end{array}
\right .
$$
The rest of the proof follows as in the classical case (see \cite{ar},p.57-58).
\end{proof}

\begin{rema}
The special form of condition (\ref{tech3}) comes from the two following constraints : one must have
$\beta \geq \alpha$ in order to preserve the regularity of perturbed curves $\gamma +h$, and
$\beta \geq 1-\alpha$ in order to ensure that the first quantity of equation (\ref{tech5}) goes to
zero when $\epsilon$ goes to zero. Note that $\alpha =1/2$ plays a special role for these sets of
conditions, as this is the only one for which the regularity of curves and variations are equals.
\end{rema}

\section{Proof of lemma \ref{simpli}}

This comes essentially from the integral formula (\ref{formint}). Indeed,
$$
\int_a^b \di {\Box_{\epsilon} \over \Box t} (f_{\epsilon} (t)h (t) ) dt
$$
is a combination of the following quantities
$$
\di {1\over 2\epsilon} \di\int_t^{t+\sigma \epsilon} f_{\epsilon} (s) h(s) ds ,\ \sigma=\pm ,
$$
for $t=a$ or $t=b$.\\

As $h(a)=h(b)=0$ and $h\in C^{\beta} (a,b)$, we have for $t=a$ or $t=b$,
$$\sup_{s\in \{ t,t+\sigma \epsilon \} } \mid h(s) \mid =
\sup_{s\in \{ t,t+\sigma \epsilon \} } \mid h(s) -h(t) \mid \leq C \epsilon^{\beta } ,$$
for some constant $C$. Moreover, using condition (\ref{tech2}), we easily obtain
$$
\sup_{s\in \{ t,t+\sigma \epsilon \} } \mid f_{\epsilon} (s) h(s) \mid C' \epsilon^{\alpha+\beta -1} .
$$
where $C'$ is a constant. Using this inequality, we deduce
$$
\mid \int_a^b \di {\Box_{\epsilon} \over \Box t} (f_{\epsilon} (s)h (s) ) ds \mid
=O(\epsilon^{\alpha +\beta -1}).
$$
We only prove the second inequality of equation (\ref{tech5}) for $Op_{\epsilon} =\Box_{\epsilon}$ and
$Op'_{\epsilon} =\Box_{\epsilon}$. The remaining cases are proved in the same way. \\

As $h\in C^{\beta} (a,b)$, we have
$$
\sup_{t\in [a,b]} \mid \Box_{\epsilon} h (t)\mid \leq C\epsilon^{\beta} ,
$$
for some constant $C$ (see lemma \ref{regsd}). Moreover, using (\ref{tech2}), we obtain
$$
\sup_{s\in [a,b] } \mid \Delta_{\epsilon}^{\sigma} (f_{\epsilon})(s) \mid \leq C^{\sigma}
\epsilon^{\alpha -1},
$$
for some constant $C^{\sigma}$, $\sigma =\pm$. We deduce
$$
\sup_{s\in [a,b]} \mid \Box_{\epsilon} f_{\epsilon} (s)\mid \leq C' \epsilon^{\alpha -1} .
$$
As a consequence, we obtain the inequality
$$
\mid \epsilon \int_a^b \di {\Box_{\epsilon} f_{\epsilon} \over \Box t} \di {\Box_{\epsilon} h \over \Box t}
dt \mid
\leq C" \epsilon^{\alpha +\beta } ,
$$
for some constant $C"$. This concludes the proof of lemma \ref{simpli}.

\section{Application : least action principle and non-linear Schr\"odinger equations}

\subsection{Least action principle and the Schr\"odinger equation}

In this section we gives a variational principle whose extremals are solutions of the {\it Schr\"odinger
equation}.\\

We consider the following non-linear Schr\"odinger's equation (obtained in \cite{bc2},\cite{cr2}) :

\begin{equation}
\label{gse}
\left [ 2i\gamma m \left [ -\di {1\over \psi} \di \left (
{\partial \psi \over \partial x}\right ) ^2 \left ( \di i\gamma
+{a_{\epsilon} (t)\over 2} \right ) +\di {\partial \psi \over
\partial t} +\di {a_{\epsilon} (t)\over 2} \di {\partial^2 \psi
\over \partial x^2} \right ] =(U(x) +\alpha (x) )\psi \right ]_{\epsilon}
,
\end{equation}
where $m>0$, $\gamma \in \rR$, $U:\rR \rightarrow \rR$, $a_{\epsilon} :\rR \rightarrow \cC$, $\alpha (x)$ is an
arbitrary continuous function.

The main result of this section is an analogue of the {\it Hamilton's principle of least action}
(see \cite{ar},p.59) for (\ref{gse}).

\begin{thm}
\label{nse}
Solutions of the non-linear Schr\"odinger equation (\ref{gse}) coincide with extremals of the
functional associated to
\begin{equation}
\label{gnew}
L(x(t),\Box_{\epsilon} x(t),t)=(1/2)m(\Box_{\epsilon} x(t))^2 +U(x),
\end{equation}
on the space of $C^{1/2}$ curves, where $x(t)$ and $\psi_{\epsilon} (x,t)$ are related by
\begin{equation}
\label{conec}
\di {\Box_{\epsilon} x \over \Box t}=-i2\gamma \di {\partial \ln (\psi (x,t))\over \partial x} ,
\end{equation}
and if $a_{\epsilon} (t)$ is such that
\begin{equation}
\label{aeps}
a_{\epsilon } (t)=\di {1\over 2} \left [
(\Delta_{\epsilon}^+ x(t))^2 - (\Delta_{\epsilon}^- x(t))^2 \right
] -i \di {1\over 2} \left [ (\Delta_{\epsilon}^+ x(t))^2 +
(\Delta_{\epsilon}^- x(t))^2 \right ] .
\end{equation}
\end{thm}

\begin{rema}
i. The nonlinear Schr\"odinger equation (\ref{gse}) was derived in \cite{bc2} using an analogue of
the Euler-Lagrange equation (\ref{stat}) proposed by Nottale \cite{no} in the context of the Scale
relativity theory. This derivation was done in the framework of the {\it local fractional calculus}
developped in \cite{bc1} and under an assumption concerning the existence of solutions to
a particular fractional differential equation. However, as proved in (\cite{cr2},part I,$\S$.4.3,
\cite{bc3}) such assumptions can't be satisfied.\\

ii. In \cite{cr2}, equation (\ref{gse}) was derived using a ``scale quantization procedure", which
gives a way to pass from classical mechanics to quantum mechanics, avoiding the problems of \cite{bc2}.
However, the Euler-Lagrange equation used in \cite{cr2} comes from scale quantization, which is
an abstract and formal way to derive the analogue of (\ref{stat}) from the classical Euler-Lagrange
equation (see $\S$.\ref{srt}).
\end{rema}

\begin{proof}
As $\partial L/\partial \Box_{\epsilon} x= m\Box_{\epsilon} x$, its differential is given by
$$
D(\partial L /\partial \Box_{\epsilon} x )=(0,m,0) ,
$$
so that condition (\ref{tech1}) is satisfied. By theorem \ref{eulerlagequa}, extremals of our
functional satisfy the Euler-Lagrange equation
\begin{equation}
\label{qele}
\left [ m\di {\Box_{\epsilon } \Box_{\epsilon} x(t)\over  \Box t} =\di {dU \over dx} (x)
\right ]_{\epsilon} .
\end{equation}
We denote
$$
f(x,t)=\di {\partial \ln (\psi (x,t))\over \partial x} (x,t) .
$$
We apply theorem \ref{tech} with $n=2$, in order to compute $\Box_{\epsilon} f(x(t),t) /\Box t$. We have
\begin{equation}
\label{cal1}
\left .
\begin{array}{lll}
\di {\Box_{\epsilon} \over \Box t} \di \left (
\di {\partial \ln (\psi ) \over \partial x} (x(t),t) \right )
& = & \di {\Box_{\epsilon} x\over \Box t} \di {\partial \over \partial x} \left (
\di {\partial \ln (\psi (x,t))\over \partial x}
\right ) (x(t),t) \\
 & & +
\di {\partial \over \partial t}
\left (
\di {\partial \ln (\psi (x,t))\over \partial x}
\right ) (x(t),t) \\
 & & +\di {1\over 2} a_{\epsilon} (t) \di {\partial^2 \over \partial x^2}
\left (
\di {\partial \ln (\psi (x,t))\over \partial x}
\right ) (x(t),t)
+o(\epsilon^{1/2} ) .
\end{array}
\right .
\end{equation}
Elementary calculus gives
$$
\di {\partial \ln (\psi (x,t))\over \partial x} =\di {1\over \psi} \di {\partial \psi\over
\partial x},\ \ \mbox{\rm and}\ \ \
\di {\partial \over \partial x}
\left (
{1\over \psi} \di {\partial \psi\over
\partial x}
\right ) =
{1\over \psi} \di {\partial^2 \psi\over
\partial^2 x}
-
{1\over \psi^2} \di \left ( {\partial \psi\over
\partial x}\right ) ^2 .
$$
Hence, we obtain
$$
\left .
\begin{array}{lll}
\di {\Box_{\epsilon} x\over \Box t} \di {\partial \over \partial x} \left (
\di {\partial \ln (\psi (x,t))\over \partial x}
\right ) (x(t),t) & = &
-i2\gamma \di {\partial \ln (\psi )\over \partial x}
\di {\partial \over \partial x}
\left (
\di {\partial \ln (\psi )\over \partial x} \right )
(x(t),t) , \\
 & = & -i\gamma \di {\partial \over \partial x}
 \left [
 \di \left (
 \di {\partial \ln (\psi )\over \partial x}
 \right )
 ^2
 \right ]
(x(t),t) , \\
 & = & -i\gamma \di {\partial \over \partial x}
 \left [
 \di {1\over \psi^2}
 \di \left (
 \di {\partial \psi \over \partial x}
\right ) ^2
\right ]
(x(t),t) .
\end{array}
\right .
$$
We then have
$$
\left .
\begin{array}{l}
\di {\Box_{\epsilon} \over \Box t}
\left (
\di {\partial \ln (\psi (x,t))\over \partial x} (x(t),t)
\right )  \\
= \di {\partial \over \partial x}
\left [
-i\gamma \di {1\over \psi^2}
\left (
\di {\partial \psi \over \partial x}
\right ) ^2
+\di {\partial \ln (\psi ) \over \partial t}
  +\di {1\over 2} a_{\epsilon} (t)
\left [
\di {1\over \psi}
\di {\partial^2 \psi \over \partial x^2} -\di {1\over \psi^2}
\left (
\di {\partial \psi \over \partial x}
\right ) ^2
\right ]
\right ]
+o(\epsilon^{1/2} ) ,\\
= \di {\partial \over \partial x} \left [
 -{1\over \psi^2}
\left ( \di {\partial \psi\over \partial x} \right ) ^2
\left (
i\gamma +\di {a_{\epsilon} (t) \over 2} \right )
+\di {1\over \psi} \di {\partial \psi\over\partial t}
+\di {a_{\epsilon} (t) \over 2} \di {1\over \psi} \di
{\partial^2 \psi \over \partial x^2}
\right ] +o(\epsilon^{1/2} ).
\end{array}
\right .
$$

As a consequence, equation (\ref{cal1}) is equivalent to
$$
\di {\partial\over \partial x}
\left [
i2\gamma m
\left [
-{\di 1\over \psi^2} \left (
\di {\partial \psi \over \partial x}
\right )
^2
\left (
i\gamma  +\di {a_{\epsilon} (t) \over 2}
\right )
+\di {1\over \psi} \di {\partial \psi\over\partial t} \right ]
+\di {a_{\epsilon} (t) \over 2} \di {1\over \psi} \di {\partial^2 \psi \over
\partial x^2}
\right ] =\di {\partial U\over \partial x}.
$$
By integrating with respect to $x$, we obtain
$$
i2\gamma m\left [
-{\di 1\over \psi^2} \left (
\di {\partial \psi \over \partial x}
\right )
^2
\left (
i\gamma  +\di {a_{\epsilon} (t) \over 2}
\right )
+\di {1\over \psi} \di {\partial \psi\over\partial t} \right ]
+\di {a_{\epsilon} (t) \over 2} \di {1\over \psi} \di {\partial^2 \psi \over
\partial x^2}
=U(x)+\alpha (x)+o(\epsilon^{1/2} ),
$$
where $\alpha (x)$ is an arbitrary function. This concludes the proof.
\end{proof}

A great deal of efforts have been made in order to generalize the classical linear Schr\"odinger
equation (see for example De Broglie \cite{debr1},\cite{debr2} and Lochak \cite{lo}). However, these
generalizations are in general ad-hoc one, choosing some particular non linear terms in order
to solve some specific problems of quantum mechanics (see for example \cite{bbm},\cite{pard},
\cite{puz}). On the contrary, the non differentiable least action principle impose a fixed non linear
term.\\

In order to recover the classical linear Schr\"odinger equation, we must specialize the functional
space on which we work. Precisely, we have :

\begin{thm}
Solutions of the Schr\"odinger equation
\begin{equation}
\left [ i\bar{h} \di {\partial \psi\over \partial t} +\di {\bar{h}^2 \over 2m} \di {\partial ^2
\psi \over \partial x^2} =U(x)\psi \right ]_{\epsilon} ,
\end{equation}
where $\bar{h}=h/2\pi$, coincide with extremals of the functional associated to
\begin{equation}
\label{gnew}
L(x(t),\Box_{\epsilon} x(t),t)=(1/2)m(\Box_{\epsilon} x(t))^2 +U(x),
\end{equation}
on the space of $C^{1/2}$ curves $\gamma :\ x=x(t)$ satisfying,
\begin{equation}
\label{regcond}
\di {1\over 2} \left [
(\Delta_{\epsilon}^+ x(t))^2 - (\Delta_{\epsilon}^- x(t))^2 \right
] -i \di {1\over 2} \left [ (\Delta_{\epsilon}^+ x(t))^2 +
(\Delta_{\epsilon}^- x(t))^2 \right ]=-i\bar{h}/m ,
\end{equation}
where $x(t)$ and $\psi_{\epsilon} (x,t)$ are related by
\begin{equation}
\label{conec}
\di {\Box_{\epsilon} x \over \Box t}=-i\di {\bar{h}\over m} \di {\partial \ln (\psi (x,t))\over \partial x} ,
\end{equation}
\end{thm}

\begin{proof}
This follows easily from the calculations made in the proof of theorem \ref{nse}.
\end{proof}

For different derivations of the Schr\"odinger equation, we refer to the work of Nelson on
{\it stochastic mechanics} (\cite{ne1},\cite{ne2}) and Feynman \cite{fey}, where he developp a
principle of least action, different from the one presented here.

 \subsection{About the scale relativity theory}
\label{srt}

This final section is informal and discuss the connexion between our non differentiable variational
principle and the scale relativity theory. In the following, we don't give a precise definition to the
word {\it fractal}. The only property which is assumed is that fractals are scale dependent objects.
We refer to \cite{fac} for more details.\\

The {\it scale relativity theory} developped by Nottale \cite{no}, gives up the assumption of the
differentiability of space-time by considering what he calls a {\it fractal space-time},
and extending the Einstein's principle of relativity to scales.\\

One of the consequences of such a theory is that there exists an infinity of geodesics\footnote{This
notion is not well defined, and we refer to \cite{no} for more details.} and that geodesics are
fractal curves. On such curves, one must developp a new differential calculus taking into account
the non differentiable character of the curve. The scale derivative introduced by Nottale is the
analogue of the scale derivative introduced in this paper.\\

The {\it scale relativity principle} can be state as follows : {\it The equations of physics keep
the same form under scale transformations} (see \cite{no}).\\

As a consequence, the scale relativity principle allows us to pass from classical mechanics
to quantum mechanic via a simple procedure : {\it one must change the classical derivative in
Newton's fundamental equation of dynamics by the scale derivative} (see \cite{no2}).\\

As Newton's equation is written via an Euler-Lagrange equation of the form
\begin{equation}
\di {d\over dt} \left [ \di {\partial L\over \partial v} \right ] =\di {\partial L\over \partial x} ,
\end{equation}
this procedure, called {\it scale quantization} in \cite{cr2}, gives a quantum analogue of
the form
\begin{equation}
\label{sele}
\di {\Box_{\epsilon} \over \Box t} \left [ \di {\partial L\over \partial v} \right ] =
\di {\partial L\over \partial x} ,
\end{equation}
where $v$ is of course a complex quantity defined by
\begin{equation}
v=\di {\Box_{\epsilon} x\over \Box t} .
\end{equation}
As a consequence, scale quantization gives an Euler-Lagrange equation similar to the one
obtained via the non differentiable variational principle introduced in this paper. The
non differentiable variational principle can be considered as an attempt to developp the
mathematical foundations of the scale relativity principle.


\begin{thebibliography}{15}
\bibitem{ar} Arnold V.I., {\it Mathematical methods of classical mechanics}, 2d edition, Graduate
Texts in Mathematics 60, Springer-Verlag, 1989.

\bibitem{bbm} Bialynicky-Birula I, Mycielsky J, Ann. Phys. 100, 62, 1976.

\bibitem{bc1} Ben Adda F, Cresson J, About non differentiable functions, Journ.
Mathematical Analysis and Applications 263, pp. 721-737, 2001.

\bibitem{bc2} Ben Adda F., Cresson J., Quantum derivatives and the Schr\"odinger equation, Chaos, solitons and fractals,
Vol. 19, no.5, 1323-1334, 2004.

\bibitem{bc3} Ben Adda F, Cresson J, Fractional differential equations and the Schr\"odinger
equation, 27.p, to appear in {\it Applied Mathematics and Computations}, 2004.

\bibitem{debr1} Broglie de L, {\it Non-linear wave mechanics}, Elsevier, Amsterdam, 1960.

\bibitem{debr2} Broglie de L, {\it Nouvelles perspectives en microphysique}, Coll. Champs
Flammarion, 1992.

\bibitem{cr1} Cresson J, Scale relativity for one dimensional non differentiable manifolds, Chaos,
 Solitons and fractals Vol. 14, no.4, pp. 553-562, 2002.

\bibitem{cr2} Cresson J., Scale calculus and the Schr\"odinger equation, Journal of Mathematical Physics, Vol. 44, No. 11,
4907-4938, 2003.

\bibitem{fac} Falconer K., {\it Fractal geometry. Mathematical foundations and applications},
John Wiley and Sons, 1990.

\bibitem{fey} Feynman R.P., The development of the space-time view of quantum electrondynamics,
Nobel lecture, December 11, 1965.

\bibitem{fh} Feynman R, Hibbs A, {\it Quantum mechanics and path integrals}, MacGraw-Hill,
1965.

\bibitem{lo} Lochak G, Annales de la fondation Louis de Broglie 22, no.1, p.1-22,
no.2, p.187-217, 1997.

\bibitem{mt} Milne-Thomson L.M., {\it The calculus of finite differences}, Chelsea Publ. Comp.,
1981.

\bibitem{ne1} Nelson E., {\it Dynamical theories of Brownian motion}, 2d edition, 2001, Princeton
University Press, 1967.

\bibitem{ne2} Nelson E., Derivation of the Schr\"odinger equation from Newtonian mechanics,
Physical Review 150 (1966).

\bibitem{no} Nottale L., {\it Fractal space-time and microphysics}, World Scientific, 1993.

\bibitem{no2} Nottale L., Scale-relativity and quantization of the universe I.  Theoritical
framework, Astron. Astrophys. 327, 867-899 (1997).

\bibitem{pard} Pardy M, To the nonlinear quantum mechanics, preprint 2002, arxiv:quant-ph/
0111105

\bibitem{puz} Puszkarz W, On the Staruszkiewicz modification of the Schr\"odinger equation,
preprint 1999, arxiv:quant-ph/9912006

\bibitem{spi} Spivak M, {\it A comprehensive introduction to differential geometry}, Publish or
Perish, Berkeley, 1979.

\bibitem{tri} Tricot C, {\it Courbes et dimension fractale}, 2d Ed., Springer, 1999.
\end{thebibliography}
\end{document}